\documentclass[12pt]{article}

\usepackage{xcolor,graphicx}
\usepackage{amssymb,amsmath}
\usepackage{comment}
\usepackage{enumerate}

\DeclareMathAlphabet{\eufrak}{U}{}{}{}
\SetMathAlphabet\eufrak{normal}{U}{euf}{m}{n}
\SetMathAlphabet\eufrak{bold}{U}{euf}{b}{n}

\usepackage[colorlinks=true, urlcolor=blue,linkcolor=blue, citecolor=blue]{hyperref}

\numberwithin{equation}{section}

\newenvironment{Proof}{\removelastskip\par\medskip
\noindent{\em Proof.} \rm}{\penalty-20\null\hfill$\square$\par\medbreak}

\allowdisplaybreaks

 \def\real{{\mathord{\mathbb R}}}
 
 \def\inte{{\mathord{\mathbb N}}}
 
 \def\qu{{\mathord{\mathbb Z}}}

 \def\real{{\mathord{{\rm I\kern-3pt R}}}}        
 
 \def\inte{{\mathord{{\rm I\kern-3pt N}}}}
 \def\sZZ{{\rm Z\kern-.45em{}Z}}

 \def\sQQ{{\kern 0.27em \vrule height1.45ex width0.03em depth0em
           \kern-0.30em \rm Q}}
 \def\qu{{\mathchoice
         {\sQQ}
         {\sQQ}
   {\kern 0.225em \vrule height1.05ex width0.025em depth0em \kern-0.25em \rm Q}
   {\kern 0.180em \vrule height0.78ex width0.020em depth0em \kern-0.20em \rm Q}
         }}
 \def\sGG{{\kern 0.27em \vrule height1.45ex width0.03em depth0em
           \kern-0.30em \rm G}}
 \def\gg{{\mathchoice
         {\sGG}
         {\sGG}
   {\kern 0.225em \vrule height1.05ex width0.025em depth0em \kern-0.25em \rm G}
   {\kern 0.180em \vrule height0.78ex width0.020em depth0em \kern-0.20em \rm G}
         }}

 \newtheorem{prop}{Proposition}[section]
 \newtheorem{lemma}[prop]{Lemma}

 \newtheorem{theorem}[prop]{Theorem}

 \def\P{{\mathord{\mathbb P}}}

\def\E{\mathop{\hbox{\rm I\kern-0.20em E}}\nolimits}

 \newcounter{hyp}
\title{\huge 
  A recursive algorithm for selling at the ultimate maximum in regime-switching models
}

\author{Yue Liu
  \\
{\normalsize School of Finance and Economics}\\
{\normalsize Jiangsu University}\\
{\normalsize Zhenjiang 212013}\\
{\normalsize P.R. China}
\and Nicolas Privault 
\\ 
\hskip-0.55cm {\normalsize School of Physical and Mathematical Sciences}
\\ 
{\normalsize Nanyang Technological University}
\\ 
{\normalsize 21 Nanyang Link} 
\\ 
{\normalsize Singapore 637371}
\\
{\normalsize           }\\
}

\usepackage{subcaption}

\begin{document}

\hyphenation{otherwise} 
\hyphenation{func-tio-nals}
\hyphenation{Privault}

\maketitle

\vspace{-1cm}
 
\begin{abstract}
 We propose a recursive algorithm for the numerical computation of
 the optimal value function
 $\inf_{t\le\tau\le T} \E \Big[\sup_{0\le s\le T } Y_s / Y_{\tau} \big| {\cal F}_t\Big]
 $
 over the stopping times $\tau$ with respect to the filtration
 of a geometric Brownian motion  $Y_t$ with Markovian
 regime switching. 
 This method allows us to determine the boundary functions of the optimal
 stopping set when no associated Volterra integral equation is available.
 It applies in particular when regime-switching drifts have mixed signs,
 in which case the boundary functions may not be monotone. 
\end{abstract} 


\baselineskip0.7cm

\section{Introduction}
\setcounter{equation}{0}
The study of optimal stopping of Brownian motion
as close as possible to its ultimate maximum has
been initiated in Graversen, Peskir and Shiryaev \cite{graversen2}.
For geometric Brownian motion, the optimal prediction problem
\begin{equation}
\label{vt}
V_t=\inf\limits_{t\le\tau\le T} \E \left[\sup\limits_{0\le s\le T }\frac{Y_s}{Y_{\tau}}\Big|{\cal F}^0_t\right]
\end{equation}
 of selling at the ultimate maximum
 over all $({\cal F}^t_s)_{s\in [t,T]}$-stopping times $\tau \in [t,T]$
 has been solved in \cite{dutoit} by Du Toit and Peskir
 when the asset price $(Y_t)_{t\in \real_+}$ is modeled
 by a geometric Brownian motion and $({\cal F}^t_s)_{s\in [t,T]}$
 is filtration generated by $(B_s-B_t)_{s\in [t,T]}$,
 see \cite{peskir} for background on optimal stopping and
 free boundary problems, and Chapter~VIII therein for ultimate
 position and maximum problems.
\\

\noindent 
 This framework has recently been 
 extended in Liu and Privault \cite{liu-privault}
 to the regime-switching model
\begin{equation}
\label{1gb}
dY_t=\mu(\beta_t)Y_tdt+\sigma (\beta_t)Y_tdB_t, \qquad 0 \le t\le T,
\end{equation}
 driven by a finite-state, observable continuous-time Markov chain
 $(\beta_t)_{t\in \real_+}$ with state space ${\cal M} := \{1,2, \ldots ,m\}$
 independent of the standard Brownian motion $(B_t)_{t\in \real_+}$
 on a filtered probability space
 $(\Omega, ({\cal F}_t)_{t\in \real_+}, \P)$,
 where $({\cal F}_t)_{t\in \real_+}$ is the filtration
 generated by $(B_t)_{t\in \real_+}$ and $(\beta_t)_{t\in \real_+}$
 and
 $\mu : {\cal M} \longrightarrow \real$, and
 $\sigma : {\cal M} \longrightarrow (0,\infty )$
 are deterministic functions.
\\
 
\noindent
 Regime-switching models were introduced by Hamilton \cite{hamilton}
 in the framework of time series, 
 in order to model the influence of external market factors.
 European options have been priced
 in continuous time regime-switching models
 by Yao, Zhang and Zhou \cite{yaozhangzhou}
 using a successive approximation algorithm.
 Optimal stopping for option pricing
 with regime switching has been dealt with in e.g.
 Guo~\cite{guoxin} 
 and Le and Wang~\cite{le-wang}.
\\

 It has been shown in particular in \cite{dutoit} that
 the boundary function $b(t)$ is nonincreasing and
 continuous in $t\in[0,T]$ and
 satisfies a Volterra integral equation of the form 
\begin{equation}
\label{ebdry}
G(t,b (t))=J(t,b (t))-\int_t^{T}K(t,r,b (t) )dr,
\end{equation}
 $0 \leq t \leq T$, with given terminal condition $b(T)$, 
where the functions $J(t,x)$ and $K(t,r,x )$ are
specified in \cite{dutoit}. 
\\
 
\noindent
 Under regime switching, the optimal boundary functions
 depend on the regime state of the system, and they may not be monotone
 if the drift coefficients $(\mu (i))_{i\in {\cal M}}$ 
 have switching signs, cf. Figures~\ref{fig2-3} and \ref{fig3} in
 Section~\ref{last}.
 Essentially, a boundary function increases when there is 
 sufficient time to switch from a state with negative drift to
 a state with positive drift and to 
 remain there until maturity, and is decreasing otherwise. 
 We refer to \cite{peskir-samee1} and
 \cite{peskir-samee2} for other optimal settings that involve
 non monotone boundary functions.
\\ 
 
 In the regime switching setting however,
 no Volterra equation such as \eqref{ebdry} 
 is available in general for boundary functions,
 cf. Section~5 and Remark~5.5 of \cite{liu-privault}. 
 In addition, the free boundary problem in the regime switching
 case consists in a system of interacting PDEs, making its direct
 solution more difficult, cf. Proposition~5.2 in \cite{liu-privault}. 
 In Buffington and Elliott~\cite{buffington} 
 a free boundary problem has been solved 
 under an ordering assumption
 on the boundary functions in the two-state case, 
 see Assumption~3.1 therein,
 however this condition may not hold in general in our setting, 
 cf. Figure~\ref{fig3} below,
 and their method is specific to American options.
\\
 
 \noindent
 In this paper we construct a recursive 
 algorithm for the numerical solution of \eqref{vt}
 in the regime-switching model \eqref{1gb},
 that includes the case where the drifts
 $(\mu (i))_{i\in {\cal M}}$ have nonconstant signs.
Our algorithm has a linear complexity $O(n)$ in the number $n$
of time steps, hence in the absence of regime switching
it also performs faster than the resolution of the Volterra equation, 
which has a quadratic complexity $O(n^2)$ due to the
evaluation of the integral in \eqref{ebdry}, cf. Section~\ref{last}. 
\\

 We start by recalling the main results of \cite{liu-privault}.
 From Lemma~2.1 in
 \cite{liu-privault},
 the optimal value function $V_t$ in \eqref{vt} can be written as
$$
V_t = V ( t, \hat{Y}_{0,t} / Y_t ,\beta_t),
$$
 where the function
 $V:[0,T] \times [1, \infty) \times {\cal M} \rightarrow \real_+$
 is given by
\begin{eqnarray}
\label{vtx}
 V(t,a,j) & := &
 \inf\limits_{t\le\tau\le T}\E\left[\frac{1}{{Y_\tau}}\max(a Y_t,\hat{Y}_{t,T})\;\Big|\;\beta_t=j\right]
\\
\label{vtj}
 & = &
 \inf_{t\le \tau \le T}
 \E \left[
 G \left( \tau,\beta_\tau,
 \frac{1}{Y_\tau}
 \max \left( a Y_t , \hat{Y}_{t,\tau} \right)
 \right)\; \Big| \; \beta_t = j
 \right],
\end{eqnarray}
 for $t\in[0,T]$, $j\in {\cal M}$, $a\geq 1$, with
 $\hat{Y}_{s,t} := \max_{r\in [s,t] } Y_r$, $0\leq s \leq t \leq T$, and
\begin{equation}
\label{g2}
 G(t ,a, j) : =
 \E\left[\max \left( a , \hat{Y}_{t,T} / Y_t \right)
 \;\Big|\;\beta_t =  j \right],
 \quad t\in[0,T], \  j \in{\cal M}.
\end{equation}
Here, the infimum
is taken over all $({\cal F}_s^t)_{s\in [t,T]}$-stopping times $\tau$, 
  where 
  ${\cal F}_s^t : = \sigma ( B_r - B_t, \ \beta_r \ : \ t \leq r \leq s )$, 
  $s\in [t,T]$. 
 From Proposition~3.1 in \cite{liu-privault},
 given $\beta_t = j\in {\cal M}$ and ${\hat{Y}_{0,t}}/{Y_t} = a \in[1,\infty)$,
 $t\in[0,T]$,
 the optimal stopping time for \eqref{vt},
 or equivalently for \eqref{vtj}, is the first hitting time
\begin{equation*}
 \tau_D(t,a,j) : =\inf\left\{r\geq t \ : \; \left(
 r,\frac{\hat{Y}_{0,r}}{Y_r},\beta_r \right)\in D\right\}
\end{equation*}
 of the stopping set
\begin{equation}\label{setd}
 D := \big\{
 (t,a,j)\in[0,T] \times [1,\infty) \times {\cal M} \ : \ V(t,a,j)=G(t,a,j)
 \big\}
\end{equation}
 by the process $( r , \hat{Y}_{0,r} / Y_r , \beta_r )_{r\in [t, T] }$.
\\

The stopping set $D$ defined in \eqref{setd}
 is closed, and its shape can be characterized as
 \begin{equation*}
      D=\left\{
 (t,y,j) \in [0,T] \times [1,\infty) \times {\cal M} 
 \ :\ y\ge b_D (t,j)\right\}
\end{equation*}
 in terms of the boundary functions $b_D(t,j)$ defined by
$$
 b_D (t,j): =\inf\{x \in [1,\infty) \ : \ (t,x,j)\in D\},
 \qquad t\in [0,T], \quad j \in {\cal M}, 
$$
cf. Proposition~3.2 of \cite{liu-privault}. 
\\

 If the condition $\mu(j)\geq 0$ is not satisfied for all $j\in{\cal M}$,
 then $t\mapsto b_D(t,j)$ may not be decreasing,
 cf. Figure~\ref{fig3} below.
 On the other hand, $\mu(j) \leq 0$ for all $j\in{\cal M}$
 leads to $b_D(t,j)=1$, $t\in [0,T]$, $j\in {\cal M}$,
 which corresponds to immediate exercise, cf.
 Proposition~5.3 in \cite{liu-privault}.
\\

 In this paper we construct a recursive algorithm for
 the numerical solution of the optimal stopping problem
 \eqref{vt}, by determining the stopping set $D$ from the
 values of $V(t,a,j)$ and $G(t,a,j)$,
 cf. Theorem~\ref{shutdown} and Lemma~\ref{th2.3} below.
 As this approach does not rely on the Volterra equation
 \eqref{ebdry}, it allows us in particular to determine the boundary
 function $b_D(t,j)$ without requiring the
 condition $\mu(j)\geq 0$ for all $j\in{\cal M}$ as in \cite{liu-privault},
 cf. for example Figure~\ref{fig3} below.
 In addition we do not rely
 on closed form expressions as in \cite{dutoit}
 as they are no longer available in the regime-switching
 setting.
\\

 Our algorithm extends the method of \cite{yaozhangzhou}
 as it applies not only to the computation of expectations,
 but also to optimal stopping problems. 
 However it differs from \cite{yaozhangzhou}, even when restricted to
 expectations $\E[\phi(Y_T)]$ of payoff functions $\phi(Y_T)$,
 where $(Y_t)_{t\in[0,T]}$ follows
 \eqref{1gb}.
 In particular, the recursion of \cite{yaozhangzhou}
 is based on the jump times of the Markov chain $(\beta_t)_{t\in \real_+}$
 whereas we apply a discretization of the time interval $[0,T]$,
 and our algorithm requires the Monte Carlo method 
 only for the estimation of \eqref{iteration1} below.
\section{Main results} 
In the sequel we let $\delta_n:=T/n$, $t^n_k:=k\delta_n$,
$k=0,1,\ldots , n$, 
 ${\cal T}_n := ( t^n_0 ,t^n_1 , \ldots , t^n_n )$,
 and 
 \begin{equation*}
    \lceil s\rceil_n : = \min \big\{ t \in {\cal T}_n \ : \ t \geq s \big\},
 \quad
 s\in [0,T], \quad n \geq 1.
\end{equation*}
 In the following Theorem~\ref{shutdown}, which is proved in Section~\ref{s2},
 the function $V_n(t,a,j)$ is computed
 by the backward induction \eqref{iteration1} starting from
 the terminal time $T$.
\begin{theorem}\label{shutdown}
\noindent
$(i)$ For all 
 $t\in[0,T]$, $j=1,2,\ldots , m$ and $a\geq 1$, the solution $V(t,a,j)$ of \eqref{vtx} satisfies
\begin{equation}\label{limit}
V(t,a,j)=\lim\limits_{n\rightarrow\infty}V_n(\lceil t\rceil_n,a,j),
\end{equation} 
 where $V_n(t^n_k,a,j)$ is the discrete infimum
\begin{equation}\label{gam3}
V_n(t^n_k,a,j) := \inf\limits_{t^n_{k+1}\le\tau_n\le T}\E\left[
\frac{\hat{Y}_{0,T}}{Y_{\tau_n}}\;\Big|\;\frac{\hat{Y}_{0,t^n_k}}{Y_{t^n_k}}=a, \ \beta_{t^n_k}=j \right],\quad k=0,1,\ldots,n-1,
\end{equation}
 taken over all
 ${\cal T}_n$-valued
 stopping times $\tau_n$, and
 $V_n(T,a,j):=V(T,a,j)=a$.
\\
\noindent
$(ii)$
 The value of $V_n ( t^n_k,a,j)$ in \eqref{gam3} can be computed by
 the backward induction
\begin{align}
\label{iteration1}
V_n\left(t^n_{k-1},a,j\right)=&\E\left[G\left(t^n_k, \frac{\hat{Y}_{0,t^n_k}}{Y_{t^n_k}},\beta_{t^n_k}\right)\wedge V_n\left(t^n_k,\frac{\hat{Y}_{0,t^n_k}}{Y_{t^n_k}},\beta_{t^n_k}\right)
\;\Big|\; \frac{\hat{Y}_{0,t^n_{k-1}}}{Y_{t^n_{k-1}}}=a, \ \beta_{t^n_{k-1}}=j \right],
\end{align}
for $k=1,2,\ldots,n$, 
 under the terminal condition $V_n(T,a,j)=G(T,a,j)=a$, $a\geq 1$,
 where $G(t,a,j)$ in defined in \eqref{g2}.
\end{theorem}
 In addition, by the following Theorem~\ref{conver} proved in Section~\ref{s4},
 we provide a way to approximate the function $G(t,a,j)$
 used in \eqref{iteration1}.
In the sequel we denote by
\begin{equation}\label{density}
\varphi_r (x,y) :=\sqrt{\frac{2}{\pi}}\frac{(2y-x)}{r^{3/2}}e^{- (2y-x)^2 / 2r },
 \qquad
 0 \leq x \leq y, \ r \in(0,T],
\end{equation}
 the joint probability density function of
 $\displaystyle \left( B_r , \sup\limits_{0\leq s\leq r}B_s \right)$,
 and we let $Q := [q_{i,j}]_{1\leq i,j \leq m}$ denote
 the infinitesimal generator of $(\beta_t)_{t\in [0,T]}$, and define
\begin{equation}\label{defu}
u(j):=\mu(j) / \sigma(j) -
 \sigma(j) / 2, \qquad j \in {\cal M}.
\end{equation}
Next, we show in Theorem \ref{conver} that $G$ is approximated by a limiting
sequence $(G_n)_{n\in \inte}$
given by the backward induction \eqref{psic0} below.
\begin{theorem}
\label{conver}
 For any $t\in [0,T]$ and $j\in {\cal M}$ we have
 \begin{equation*}
      G(t,a,j)= \lim\limits_{n\rightarrow\infty}
   G_n(\lceil t\rceil_n,a,j),
\end{equation*} 
 where the limit is uniform in $a\geq 1$ 
 and $G_n(t^n_{k},a,j)$ is defined by the backward induction
\begin{align} 
\label{psic0}
&  G_n(t^n_{k-1},a,j) =
\\
\nonumber
&
e^{q_{j,j}\delta_n}
 \int_0^{\infty}\int_{-\infty}^ye^{ ( u(j)+
 \sigma(j)) x- u^2(j)T / ( 2n ) }
 G_n\left(
t^n_{k} , a\vee(\sigma(j)y)-\sigma(j)x,j \right)
 \varphi_{\delta_n}(x,y)dx dy
\\
\nonumber
 & +
 \sum_{i=1\atop i \not= j}^m
 q_{j,i}
 \int_0^{\delta_n} e^{q_{j,j}r} 
 \int_0^{\infty}\int_{-\infty}^ye^{ ( u(j)+\sigma(j) ) x- u^2(j)r / 2 }
 G_n\left(t^n_{k},a\vee(\sigma(j)y)-\sigma(j)x , i\right)
 \varphi_r (x,y)dx dy dr,
\end{align} 
 $k=1,2,\ldots,n$, 
with the terminal condition $G_n(T,a,j) =a$, $j\in {\cal M}$, $a\geq 1$.
\end{theorem}
In the particular case of constant drift $\mu$ and volatility $\sigma$ cf,
 Theorems~\ref{shutdown} and \ref{conver} also
 provide an alternative numerical solution 
 in the geometric Brownian motion model of \cite{dutoit}. 
 In this case, $(V_n(t^n_{k-1},a))_{k=1,2,\ldots,n}$ is
 computed
 from \eqref{iteration1} by the backward induction
\begin{eqnarray*}
\nonumber
 V_n\left(t^n_{k-1},a\right)
 = \int_0^\infty \int_{-\infty}^y
 \! \!
 G\left(t^n_k, e^{\sigma(\frac{\log a}{\sigma}\vee y-x)}\right)\wedge V_n\left(t^n_k,e^{\sigma(\frac{\log a}{\sigma}\vee y-x)}\right)e^{\lambda x-
 \lambda^2\delta_n / 2 } \varphi_{\delta_n} (x,y) dx dy,
\\
\end{eqnarray*}
\vskip-0.6cm
 with
\begin{align}\label{newpsi}
\quad G(t,a) =\E\left[a\vee e^{\sigma S^\lambda_{T-t}}\right]=\int_0^\infty \int_{-\infty}^y e^{ ( \log a ) \vee(\sigma y)+\lambda x- \lambda^2(T-t) / 2 }
 \varphi_{T-t}(x,y)dxdy,
\end{align}
for all $t\in[0,T]$ and $a\geq 1$, where
$S^\lambda_t:=\max\limits_{0\leq s\leq t} ( B_s +\lambda s )$, 
 $\lambda:=\mu/\sigma-\sigma/2$, and
 $\varphi_r (x,y)$ is given by \eqref{density}.
In the general regime switching setting,
the function $G(t,a)$ in \eqref{newpsi} 
can be estimated by Monte Carlo,
while in the absence of regime switching 
it can be computed in closed form, cf. (2.7) in \cite{dutoit}.
\\
 
 In Sections~\ref{s2} and \ref{s4}
 we prove Theorems~\ref{shutdown} and \ref{conver}.
 Numerical illustrations are presented in Section~\ref{last}
 with and without regime switching.
 We observe in particular that boundary functions may not be
 monotone when the drift coefficients
 $\mu(j)$, $j\in {\cal M}$, have different signs. 
\section{Proof of Theorem~\ref{shutdown}}
\label{s2}
\setcounter{equation}{0}
\noindent
$(i)$ First, we note that
for any $({\cal F}^t_s)_{s\in [t,T]}$-stopping time $\tau \in [ t ,T]$ we have
\begin{eqnarray}
\nonumber
\E\left[\frac{(aY_{ t })\vee\hat{Y}_{{ t },T}}{Y_{\tau}}\;\Big|\;\beta_{ t }=j\right]
 & = & \E\left[\frac{a\vee\left(
 \hat{Y}_{{ t },T} / Y_{ t }
 \right)}{
 Y_{\tau} / Y_{ t } }\;\Big|\; \frac{\hat{Y}_{0,{ t }}}{Y_{ t }}=a, \ \beta_{ t }=j\right]
\\
\label{condt}
& = & \E\left[
\frac{\hat{Y}_{0,T}}{Y_{\tau}}\;\Big|\; \frac{\hat{Y}_{0,{ t }}}{Y_{ t }}=a, \ \beta_{t}=j \right]
,
\end{eqnarray}
 $t\in[0,T]$, $j\in {\cal M}$, $a\in [1,\infty)$,
 since ${\hat{Y}_{0,t }}/{Y_{ t }}$ is conditionally independent of
$$
 \left(\frac{Y_{ t }}{Y_\tau},\frac{\hat{Y}_{{ t },T}}{Y_\tau}\right)
 = \left(\exp\left(-\int_t^\tau \sigma(\beta_r )d\tilde{B}_r \right),\exp\left(
 \sup_{t \leq v \leq T} \int_t^v \sigma(\beta_r )d\tilde{B}_r
 -\int_t^\tau \sigma(\beta_t )d\tilde{B}_r \right)\right)
$$
 given $\beta_{ t }$, where $(\tilde{B}_v)_{v\in [0,T]}$ is the
 drifted Brownian motion 
 \begin{equation}
   \label{dfkjlsdf}
   \tilde{B}_v : = B_v + \int_0^v u(\beta_r)dr,\qquad v \in [0,T],
   \end{equation}
 and $u(j):={\mu(j)}/{\sigma(j)}-{\sigma(j)}/{2}$,
 $j\in {\cal M}$, is defined in \eqref{defu}.
 Hence by \eqref{vtx} we have
\begin{eqnarray*}
V(t^n_k,a,j)&=&\inf\limits_{t^n_k\le\tau\le T}\E\left[
\frac{\hat{Y}_{0,T}}{Y_\tau }\;\Big|\; \frac{\hat{Y}_{0,{t^n_k}}}{Y_{t^n_k}}=a, \ \beta_{t^n_k}=j \right]\nonumber\\
&\leq&\inf\limits_{{t^n_k}\le\tau_n\le T}\E\left[
\frac{\hat{Y}_{0,T}}{Y_{\tau_n} }\;\Big|\;\frac{\hat{Y}_{0,{t^n_k}}}{Y_{t^n_k}}=a, \ \beta_{t^n_k}=j\right]\nonumber\\
&\leq&\inf\limits_{{t^n_{k+1}}\le\tau_n\le T}\E\left[
\frac{\hat{Y}_{0,T}}{Y_{\tau_n} }\;\Big|\; \frac{\hat{Y}_{0,{t^n_k}}}{Y_{t^n_k}}=a, \ \beta_{t^n_k}=j \right]\nonumber\\
&=&V_n({t^n_k},a,j),
\end{eqnarray*}
 $k=0,1,\ldots,n-1$, $j \in {\cal M}$, $a\geq 1$,
 where we used \eqref{gam3} and the infimum is taken over all
 ${\cal T}_n$-valued discrete stopping times $\tau_n$.
 Therefore by the continuity of $V(t,a,j)$ with respect to $t$,
 cf. e.g. \cite{peskir}, Chap III, \S 7.1.1 page~130 and \S 7.4.1 pages~135-136, we obtain
\begin{equation}
\label{3<}
 V(t,a,j) =
 \lim\limits_{n\rightarrow\infty} V(\lceil t\rceil_n,a,j)
 \leq
 \liminf\limits_{n\rightarrow\infty}V_n(\lceil t\rceil_n,a,j).
\end{equation}
\noindent
$(ii)$ On the other hand, by \eqref{condt} we have
\begin{eqnarray}
\nonumber
\lefteqn{
\limsup\limits_{n\rightarrow\infty}V_n({\lceil t\rceil_n},a,j)
 = \limsup\limits_{n\rightarrow\infty}\inf\limits_{{\lceil t\rceil_n}+\delta_n\le\tau_n\le T}\E\left[
\frac{\hat{Y}_{0,T}}{Y_{\tau_n}}\;\Big|\; \frac{\hat{Y}_{0,{\lceil t\rceil_n}}}{Y_{\lceil t\rceil_n}}=a, \ \beta_{\lceil t\rceil_n}=j \right]
}
\\
\nonumber
&=& \limsup\limits_{n\rightarrow\infty}\inf\limits_{{\lceil t\rceil_n}+\delta_n\le\tau_n\le T}\E\left[\frac{(aY_{\lceil t\rceil_n})\vee\hat{Y}_{{\lceil t\rceil_n},T}}{Y_{\tau_n}}\;\Big|\;\beta_{\lceil t\rceil_n}=j\right]
\\
\nonumber
 & = &
 \limsup\limits_{n\rightarrow\infty}
 \sum\limits_{l=1}^m
 \left[e^{(\lceil t\rceil_n-t)Q}\right]_{j,l}\inf\limits_{{\lceil t\rceil_n}+\delta_n\le\tau_n\le T}\E\left[\frac{(aY_{\lceil t\rceil_n})\vee\hat{Y}_{{\lceil t\rceil_n},T}}{Y_{\tau_n}}\;\Big|\;\beta_{\lceil t\rceil_n} =l\right]
\\
\nonumber
& \leq &
 \limsup\limits_{n\rightarrow\infty}
 \inf\limits_{{\lceil t\rceil_n}+\delta_n\le\tau_n\le T}\sum\limits_{l=1}^m
 \left[e^{(\lceil t\rceil_n - t )Q}\right]_{j,l}\E\left[\frac{(aY_{\lceil t\rceil_n})\vee\hat{Y}_{{\lceil t\rceil_n},T}}{Y_{\tau_n}}\;\Big|\;\beta_{\lceil t\rceil_n} =l\right]
\\
\label{djklsdfds}
& = &
 \limsup\limits_{n\rightarrow\infty}
 \inf\limits_{{\lceil t\rceil_n}+\delta_n\le\tau_n\le T}\E\left[\frac{(aY_{\lceil t\rceil_n})\vee\hat{Y}_{{\lceil t\rceil_n},T}}{Y_{\tau_n}}\;\Big|\;\beta_t=j\right],
\end{eqnarray}
 $t\in[0,T-\delta_n ]$, $j\in {\cal M}$, $a\in [1,\infty)$,
 where $Q = [q_{i,j}]_{1\leq i,j \leq m}$ is the infinitesimal
 generator of $(\beta_t)_{t\in [0,T]}$.
 Next, we note that for every stopping time $\tau\in[t,T]$
 we have $|\lceil\tau\rceil_n-\tau|<1/n$,
 hence $(\lceil\tau\rceil_n)_{n\geq 1}$
 converges to $\tau$ uniformly in $L^\infty(\Omega)$ and pointwise.
 Hence
 we have
\begin{equation}\label{lebsg}
\lim\limits_{n\rightarrow\infty}\E\left[\frac{(aY_{\lceil t\rceil_n})\vee\hat{Y}_{{\lceil t\rceil_n},T}}{Y_{\lceil\tau\vee(t+\delta_n)\rceil_n}}\;\Big|\;\beta_{t}=j\right]
=
\E\left[\lim\limits_{n\rightarrow\infty}\frac{(aY_{\lceil t\rceil_n})\vee\hat{Y}_{{\lceil t\rceil_n},T}}{Y_{\lceil\tau\vee(t+\delta_n)\rceil_n}}\;\Big|\;\beta_{t}=j\right]
,
\end{equation}
$t\in [0,T-\delta_n]$,
for any stopping time $\tau\in[t,T]$, where we applied
 Lebesgue's dominated convergence theorem based
 on the bounds \eqref{lebg2} and \eqref{lebg3}
 stated at the end of this section.
 Hence from \eqref{djklsdfds} and \eqref{lebsg} we find,
 for any stopping time $\tau\in[t,T]$,
\begin{eqnarray*}
\nonumber
\limsup\limits_{n\rightarrow\infty}V_n({\lceil t\rceil_n},a,j)
 & \leq &
 \limsup\limits_{n\rightarrow\infty}
 \inf\limits_{{\lceil t\rceil_n}+\delta_n\le\tau_n\le T}\E\left[\frac{(aY_{\lceil t\rceil_n})\vee\hat{Y}_{{\lceil t\rceil_n},T}}{Y_{\tau_n}}\;\Big|\;\beta_t=j\right]
\\
\nonumber
& \leq &
\lim\limits_{n\rightarrow\infty}\E\left[\frac{(aY_{\lceil t\rceil_n})\vee\hat{Y}_{{\lceil t\rceil_n},T}}{Y_{\lceil\tau\vee(t+\delta_n)\rceil_n}}\;\Big|\;\beta_{t}=j\right]
\\
 & = & \E\left[\lim\limits_{n\rightarrow\infty}\frac{(aY_{\lceil t\rceil_n})\vee\hat{Y}_{{\lceil t\rceil_n},T}}{Y_{\lceil\tau\vee(t+\delta_n)\rceil_n}}\;\Big|\;\beta_{t}=j\right]\\
& = & \E\left[\frac{(aY_{t})\vee\hat{Y}_{t,T}}{Y_{\tau}}\;\Big|\;\beta_{t}=j\right]\\
&= & \E\left[\frac{\hat{Y}_{0,T}}{Y_\tau}\;\Big|\; \frac{\hat{Y}_{0,t}}{Y_t}=a, \ \beta_t=j \right],
\end{eqnarray*}
 where we applied \eqref{condt} and
 the pathwise continuity of $(Y_t)_{t\in [0,T]}$.
 Hence by \eqref{gam3}, we obtain
\begin{eqnarray*}
\limsup\limits_{n\rightarrow\infty}V_n({\lceil t\rceil_n},a,j)
\leq  \inf\limits_{t\leq\tau\leq T} \E\left[\frac{\hat{Y}_{0,T}}{Y_\tau}\;\Big|\; \frac{\hat{Y}_{0,t}}{Y_t} =a, \ \beta_t=j \right]
=V(t,a,j),
\end{eqnarray*}
$t\in [0,T-\delta_n]$,
which completes the proof of \eqref{limit} by \eqref{3<}.
\\

\noindent
$(iii)$ In order to prove \eqref{iteration1} for 
 $0 \leq k \leq l \leq n$, 
 we consider an optimal stopping time
 $\tau_n^{(t^n_l)}$ such that
\begin{equation}
\label{fjoisfd}
\E\left[\frac{\hat{Y}_{0,T}}{Y_{\tau_n^{(t^n_l)}}}\;\Big|\; \frac{\hat{Y}_{0,{t^n_k}}}{Y_{t^n_k}}=a , \ \beta_{t^n_k}=j \right]=\inf\limits_{t^n_l\le\tau_n\le T}\E\left[
\frac{\hat{Y}_{0,T}}{Y_{\tau_n}}\;\Big|\; \frac{\hat{Y}_{0,{t^n_k}}}{Y_{t^n_k}}=a, \ \beta_{t^n_k}=j \right],
\end{equation}
 where the infimum is taken over
 the discrete ${\cal T}_n$-valued
 stopping times $\tau_n$,
 and the existence of $\tau_n^{(h)}$
 is guaranteed by Corollary~2.9 of \cite{peskir} as
 in Proposition~3.1 of \cite{liu-privault}.
 We note the induction
\begin{align}
\E\left[
\frac{\hat{Y}_{0,T}}{Y_{\tau_n^{(t^n_k)}}}\;\Big|\; \frac{\hat{Y}_{0,t^n_k }}{Y_{t^n_k}}, \ \beta_{t^n_k} \right]&=\E\left[
\frac{\hat{Y}_{0,T}}{Y_{t^n_k}}\;\Big|\; \frac{\hat{Y}_{0,{t^n_k} }}{Y_{t^n_k}}, \ \beta_{t^n_k} \right]
\wedge \E\left[
\frac{\hat{Y}_{0,T}}{Y_{\tau_n^{( t^n_{k+1} ) }}}\;\Big|\; \frac{\hat{Y}_{0,{t^n_k} }}{Y_{t^n_k}}, \ \beta_{t^n_k} \right]\nonumber\\ 
&=G\left(t^n_k,\frac{\hat{Y}_{0,{t^n_k}}}{Y_{t^n_k}},\beta_{t^n_k}\right)\wedge V_n\left({t^n_k},\frac{\hat{Y}_{0,{t^n_k} }}{Y_{t^n_k}},\beta_{t^n_k}\right),\label{split}
\end{align}
 $k=0,1,\ldots,n-1$, $a\geq 1$, where $V_n$ and $G$
 are defined in \eqref{gam3} and \eqref{g2} respectively.
 By \eqref{fjoisfd}, this yields
\begin{eqnarray*}
\lefteqn{
V_n\left(t^n_{k-1},a,j\right) = \E\left[
\frac{\hat{Y}_{0,T}}{Y_{\tau_n^{(t^n_{k})}}}\;\Big|\;
\frac{\hat{Y}_{0,t^n_{k-1}}}{Y_{t^n_{k-1}}}=a, \ \beta_{t^n_{k-1}}=j \right]
}\\
&=&\E\left[\E\left[
\frac{\hat{Y}_{0,T}}{Y_{\tau_n^{(t^n_{k})}}}\;\Big|\;\frac{\hat{Y}_{0,t^n_{k}}}{Y_{t^n_{k}}},
\frac{\hat{Y}_{0,t^n_{k-1}}}{Y_{t^n_{k-1}}}=a , \ \beta_{t^n_{k}}, \ \beta_{t^n_{k-1}}=j \right]
\;\Big|\; \frac{\hat{Y}_{0,t^n_{k-1}}}{Y_{t^n_{k-1}}}=a, \ \beta_{t^n_{k-1}}=j \right]\\
&=&\E\left[\E\left[
\frac{\hat{Y}_{0,T}}{Y_{\tau_n^{(t^n_{k})}}}\;\Big|\; \frac{\hat{Y}_{0,t^n_{k}}}{Y_{t^n_{k}}}, \ \beta_{t^n_{k}} \right]
\;\Big|\; \frac{\hat{Y}_{0,t^n_{k-1}}}{Y_{t^n_{k-1}}}=a, \ \beta_{t^n_{k-1}}=j \right]\\
&=&\E\left[G\left(t^n_{k}, \frac{\hat{Y}_{0,t^n_{k}}}{Y_{t^n_{k}}},\ \beta_{t^n_{k}}\right)\wedge V_n\left(t^n_{k},\frac{\hat{Y}_{0,t^n_{k}}}{Y_{t^n_{k}}},\ \beta_{t^n_{k}}\right)
\;\Big|\; \frac{\hat{Y}_{0,t^n_{k-1} }}{Y_{t^n_{k-1}}}=a , \ \beta_{t^n_{k-1}}=j \right],
\end{eqnarray*}
 $k=1,2,\ldots,n$, where we applied
 \eqref{split}, the Markov property of
 $( \hat{Y}_{0,t} / Y_{t} , \beta_{t} )_{t\in [0,T]}$ and the relation $V_n(T,a,j)=G(T,a,j)$.
 \hfill $\Box$
\\

\noindent
We close this section with the proof of the two bounds used for 
\eqref{lebsg} above.
\begin{enumerate}[(a)]
\item Letting $\check{Y}_{t,T}:=\min\limits_{t \le v \le T}Y_v$,
  we check that, for any stopping time $\tau$ and $a\geq 1$, we have
  the bound 
\begin{equation}\label{lebg2}
\max\left(\frac{aY_{\lceil t\rceil_n}}{Y_{\lceil\tau\vee(t+\delta_n)\rceil_n}},\frac{\hat{Y}_{\lceil t\rceil_n,T}}{Y_{\lceil\tau\vee(t+\delta_n)\rceil_n}}\right)\leq \frac{a\hat{Y}_{t,T}}{Y_{\lceil\tau\vee(t+\delta_n)\rceil_n}}\leq \frac{a\hat{Y}_{t,T}}{\check{Y}_{t,T}},
\end{equation}
 in which the right hand side is integrable for all $t\in [0,T-\delta_n ]$. 
\item
 On the other hand we have 
 $\E\left[ {\hat{Y}_{t,T}} / {\check{Y}_{t,T}}\;\Big|\;\beta_t=j\right]<\infty$
 since, using the drifted Brownian motion $(\tilde{B}_v)_{v\in [0,T]}$
 defined in \eqref{dfkjlsdf} we have, using the Cauchy-Schwarz
 inequality,
\begin{eqnarray}
\label{lebg3}
\lefteqn{
 \! \! \! \! \! \! \! \! \!
 \E\left[ \frac{\hat{Y}_{t,T}}  {\check{Y}_{t,T}}\;\Big|\;\beta_t=j\right]=\E\left[e^{\sup\limits_{t\leq r\leq T}\int_t^r \sigma(\beta_v )d\tilde{B}_v - \inf\limits_{t\leq r\leq T}\int_t^r\sigma(\beta_v )d\tilde{B}_v}\;\Big|\;\beta_t=j\right]}
\\
\nonumber
&\leq &
 \sqrt{
 \E\left[e^{2\sup\limits_{t\leq r\leq T}\int_t^r \sigma(\beta_v )d\tilde{B}_v}\;\Big|\;\beta_t=j\right]
\E\left[e^{-2\inf\limits_{t\leq r\leq T}\int_t^r \sigma(\beta_v )d\tilde{B}_v}\;\Big|\;\beta_t=j\right]}
\\
\nonumber
&\leq &
 \sqrt{\E\left[e^{2\sup\limits_{t\leq r\leq T}\int_t^r \sigma(\beta_v )d\tilde{B}_v}\;\Big|\;\beta_t=j\right]}
 \\
 \nonumber
 & < & \infty, 
\end{eqnarray} 
where we conclude to finiteness by conditioning and use of
the density \eqref{density}. 
\end{enumerate}
\section{Proof of Theorem~\ref{conver}}
\label{s4}
 We start with two lemmas.
\begin{lemma}
\label{th2.3}
 For all $k = 1,2,\ldots,n$, $j\in {\cal M}$ and $a\geq 1$,
 we have
\begin{align} 
\label{psic}
& 
 G(t^n_{k-1},a,j) =
 \\
 \nonumber
 &
 e^{q_{j,j}\delta_n}
 \int_0^{\infty}\int_{-\infty}^ye^{ ( u(j)+
 \sigma(j)) x- u^2(j)T / ( 2n ) }
 G\left(
t^n_k ,
 a\vee(\sigma(j)y)-\sigma(j)x,j
 \right)
 \varphi_{\delta_n}(x,y)dx dy
\\
\nonumber
 & +
 \sum_{i=1\atop i \not= j}^m
 q_{j,i}
 \int_0^{\delta_n}
 e^{q_{j,j}r}
 \int_0^{\infty}\int_{-\infty}^ye^{ ( u(j)+\sigma(j) ) x- u^2(j)r / 2 }
 G\left(
t^n_{k-1}+r,a\vee(\sigma(j)y)-\sigma(j)x , i\right)
 \varphi_r (x,y)dx dy dr.
\end{align} 
\end{lemma}
\begin{Proof}
 Let $\tilde{\P}$ denote the probability measure defined by
$$
 \frac{d\tilde{\P}}{d{\P}} :=
 \exp\left(-\int_0^T u(\beta_r )dB_r-\frac{1}{2}\int_0^T u^2(\beta_r) dr\right),
$$
 where $u(j)$, $j\in {\cal M}$, is defined in \eqref{defu}, 
 and $(\tilde{B}_r)_{r\in [0,T]}$
 is the standard Brownian motion under $\tilde{\P}$
 defined in \eqref{dfkjlsdf}.
 From the definition \eqref{g2} of $G(t,a,j)$ we have
\begin{eqnarray}
G(t,a,j)
&= & \E\left[a\vee\exp
 \left(
 \sup\limits_{t\le s\le T}\int_t^s\sigma(\beta_r)d\tilde{B}_r \right)
 \;\Bigg|\;\beta_t=j\right]\label{forlateuse}\\
&= & \nonumber\tilde{\E} \left[e^{ ( \log a ) \vee\sup\limits_{t\le s\le T}\int_t^s\sigma(\beta_r)d\tilde{B}_r+\int_t^Tu(\beta_r)d B_r+\frac{1}{2}\int_t^Tu^2(\beta_r)dr}\;\Big|\;\beta_t=j\right]\\
&= & \nonumber\tilde{\E} \left[e^{ ( \log a ) \vee\sup\limits_{t\le s\le T}\int_t^s\sigma(\beta_r)d\tilde{B}_r+\int_t^Tu(\beta_r)d\tilde{B}_r-\frac{1}{2}\int_t^Tu^2(\beta_r)dr}\;\Big|\;\beta_t=j\right]\\
&= & \E\left[e^{ ( \log a ) \vee\sup\limits_{t\le s\le T}\int_t^s\sigma(\beta_r)dB_r+\int_t^Tu(\beta_r)dB_r-\frac{1}{2}\int_t^Tu^2(\beta_r)dr}\;\Big|\;
\beta_t=j\right],\label{g3}
\end{eqnarray}
 which allows us to remove the drift component in the
 supremum $\sup\limits_{t\le s\le T}\int_t^s\sigma(\beta_r)dB_r$.
 Next, using \eqref{g3} we write 
\begin{equation}
\label{djklsadf}
G(t^n_k,a,j) = \Phi_n(t^n_k,a,j) + \Upsilon_n(t^n_k,a,j),
\qquad
j\in {\cal M}, \quad a\geq 1,
\end{equation}
 where
\begin{equation}\label{**1}
\Phi_n(t^n_k,a,j):=\E \left[e^{ ( \log a ) \vee\sup\limits_{t^n_k \le s\le T}\int_{t^n_k}^s\sigma(\beta_r)dB_r+\int_{t^n_k}^T u(\beta_r)dB_r-\frac{1}{2}\int_{t^n_k}^Tu^2(\beta_r)dr}
 \hskip-1.2cm
 {\bf 1}_{\{T_1(t^n_k)> t^n_{k+1} \}}\Big|\beta_{t^n_k}=j\right],
\end{equation}
with $T_1(t):=\inf\{s\geq t\,:\, \beta_s\neq\beta_t\}$ for any $t\in\real_+$,
and
\begin{equation}\label{**2}
\Upsilon_n(t^n_k,a,j):=\E \left[e^{ ( \log a ) \vee\sup\limits_{{t^n_k}\le s\le T}\int_{t^n_k}^s\sigma(\beta_r)dB_r+\int_{t^n_k}^Tu(\beta_r)dB_r-\frac{1}{2}\int_{t^n_k}^Tu^2(\beta_r)dr}
 \hskip-1.2cm
 {\bf 1}_{\{T_1(t^n_k)\leq t^n_{k+1} \}}\Big|\beta_{t^n_k}=j\right].
\end{equation}
By \eqref{**1} we have, for any $k=0,1,\ldots,n-1$,
$j\in {\cal M}$, and $a\geq 1$,
\begin{align}
 & \label{k=l1}
  \nonumber
  \Phi_n(t^n_k,a,j) 
  =
  \nonumber
 e^{q_{j,j}\delta_n}
 \int_0^{\infty}\int_{-\infty}^y
 \\
 \nonumber
 & 
 \E\left[
 e^{
 ( \log a ) \vee \left( \sigma(j)y\vee\left(\sigma(j)x
 +
 \! \! \!
 \sup\limits_{t^n_{k+1} \le s\le T}\int_{t^n_{k+1}}^s
 \sigma(\beta_r)dB_r \right) \right)
 +\int_{t^n_{k+1}}^T
 u(\beta_r)dB_r+u(j)x-\frac{1}{2}\int_{t^n_{k+1} }^Tu^2(\beta_r)dr
- u^2(j) \delta_n /2 } \right]
 \\
 \nonumber
& \times \varphi_{\delta_n}(x,y)dxdy
 \\
 \nonumber
&= 
 e^{q_{j,j}\delta_n}
 \int_0^{\infty}\int_{-\infty}^y
 \\
 \nonumber
 & 
 \E\left[
 e^{ \! \!
 (
 ( ( \log a ) \vee(\sigma(j)y)-\sigma(j)x ) \vee
 \! \! \! \! \!
 \sup\limits_{t^n_{k+1}\le s\le T}\int_{t^n_{k+1} }^s\sigma(\beta_r)dB_r
 +\sigma(j)x+\int_{t^n_{k+1}}^T
u(\beta_r)dB_r+u(j)x-\frac{1}{2}\int_{t^n_{k+1} }^Tu^2(\beta_r)dr- u^2(j) \delta_n /2 ) }
 \right]
 \\
 \nonumber
 &
 \times \varphi_{\delta_n}(x,y)dxdy\\
 & 
 = e^{q_{j,j}\delta_n}
 \int_0^{\infty}\int_{-\infty}^ye^{ ( u(j)+\sigma(j) ) x-
 u^2(j)T / ( 2n ) }
 G( t^n_{k+1},e^{ ( \log a ) \vee(\sigma(j)y)-\sigma(j)x},j)
\varphi_{\delta_n}(x,y)dxdy,
\end{align}
 where in the first equality we used the fact that the time to the first
 jump of $(\beta_s)_{s\in [t,\infty )}$ after $t$ is
 exponentially distributed with parameter $-q_{j,j} > 0$
 given $\beta_t = j$, cf. e.g. \S~10.4 in \cite{privaultbkm1}.
 Next, for $\Upsilon_n(t^n_k,a,j)$, $k=0,1,\ldots,n-1$,
 $j\in {\cal M}$, and $a\geq 1$, by \eqref{**2} we see that
\vskip-1.cm
\begin{align} 
\label{k=l}
& \Upsilon_n(t^n_k,a,j) 
= 
 \sum_{i=1\atop i \not= j}^m
 q_{j,i}
 \int_0^{\delta_n} e^{q_{j,j}r} \int_0^{\infty}\int_{-\infty}^y
\\
\nonumber
&
 \E\left[
 e^{ \! \!
 ( \log a ) \vee ( \sigma(j)y)\vee
 (
 \sigma(j)x+
 \! \! \! \!
 \sup\limits_{t^n_k +r\le s\le T}\int_{t^n_k +r}^s\sigma(\beta_z )dB_z
+\int_{t^n_k +r}^T
 u(\beta_z)dB_z +u(j)x-\frac{1}{2}\int_{t^n_k +r}^Tu^2(\beta_z )dz
 -u^2(j)r/2
 ) }
 \Big| \beta_{t^n_k+r}=i\right]
 \\
 \nonumber
& 
 \varphi_r (x,y)dxdy dr 
 \\
 \nonumber
&=
 \sum_{i=1\atop i \not= j}^m
 q_{j,i}
 \int_0^{\delta_n} e^{q_{j,j}r} \int_0^{\infty}\int_{-\infty}^y
e^{ ( u(j)+\sigma(j) ) x-{u^2(j)r}/{2}}G( t^n_k +r,e^{ ( \log a ) \vee(\sigma(j)y)-\sigma(j)x},i)
\varphi_r (x,y)dxdy dr,
\end{align} 
 where we used the conditional probability distribution
$$
 {\P} ( T_1 \in dt, \ \beta_{T_1} = i \mid \beta_0 = j)
 = {\bf 1}_{[0,\infty )}(t) q_{j,i} e^{q_{j,j} t} dt,
 \qquad i\not= j\in {\cal M},
$$
 computed from the
 exponential distribution with parameter $-q_{i,i}$
 of the first jump time $T_1$ of the
 Markov chain $(\beta_t)_{t\in \real_+}$ started at $i\in {\cal M}$
 and the transition matrix
 $( - q_{i,j} {\bf 1}_{\{ i \not= j\}} / q_{i,i} )_{i,j\in {\cal M}}$
 of the embedded Markov chain, cf. e.g. \S~10.7 of \cite{privaultbkm1}. Hence
we conclude to \eqref{psic} by \eqref{djklsadf}, \eqref{k=l1} and \eqref{k=l}.
\end{Proof}
\begin{lemma}\label{gcontin}
  For any $j\in{\cal M}$, the function
  $t\mapsto G(t,a,j)$ is uniformly continuous in $t\in [0,T]$,
  uniformly in $a\geq 1$, i.e. 
  \begin{equation}
   \label{psilim}
    \lim_{\varepsilon \to 0}
    \sup_{|t-s|\leq \varepsilon}
    \sup_{a\geq 1}
    | G(t,a,j) - G(s,a,j) | = 0,
    \qquad
    j \in {\cal M}.
  \end{equation}
\end{lemma}
\begin{Proof}
  By \eqref{forlateuse}, for all $a\geq 1$ we have 
\begin{align}\label{gx}
&|G(t,a,j)-G(s,a,j)|\nonumber\\
&=\nonumber\Bigg| \E\left[a\vee\exp\left(
 \sup\limits_{t\le v \le T}\int_t^v \sigma(\beta_r)d\tilde{B}_r \right)
 \;\Bigg|\;\beta_t=j\right] -\E\left[a\vee\exp\left(
 \sup\limits_{s\le v \le T}\int_s^v \sigma(\beta_r)d\tilde{B}_r \right)
 \;\Bigg|\;\beta_s=j\right]\Bigg|\\
&\leq \Bigg|\E\left[\exp\left(
 \sup\limits_{t\le v \le T}\int_t^v \sigma(\beta_r)d\tilde{B}_r \right)
 \;\Bigg|\;\beta_t=j\right]
  -\E\left[\exp\left(
 \sup\limits_{s\le v \le T}\int_s^v \sigma(\beta_r)d\tilde{B}_r \right)
 \;\Bigg|\;\beta_s=j\right]\Bigg|, 
\end{align}
hence it suffices to show the continuity in $t\in [0,T]$
of the above bound. Similarly to \eqref{g3}, we have 
\begin{equation*}
  \E\left[\exp\left(
 \sup\limits_{t\le v\le T}\int_t^v\sigma(\beta_r)d\tilde{B}_r \right)
 \;\Bigg|\;\beta_t=j\right]
=\E\left[e^{\sup\limits_{t\le v\le T}\int_t^v\sigma(\beta_r)dB_r+\int_t^Tu(\beta_r)dB_r-\frac{1}{2}\int_t^Tu^2(\beta_r)dr}\;\Big|\;
\beta_t=j\right],
\end{equation*}
 $t\in[0,T]$, $j\in{\cal M}$.
Next, for any $n\geq 1$ we have
\begin{eqnarray*}
 \sup\limits_{\lceil t\rceil_n \le s\le T}\int_{\lceil t\rceil_n }^s\sigma(\beta_r)dB_r
 & = & \sup\limits_{\lceil t\rceil_n \le s\le T}\int_t^s\sigma(\beta_r)dB_r-\int_t^{\lceil t\rceil_n}\sigma(\beta_r)dB_r
\\
 & \leq &
 \sup\limits_{t \le s\le T}\int_{t}^s\sigma(\beta_r)dB_r
 -
 \inf\limits_{t \le s\le T}\int_{t}^{s}\sigma(\beta_r)dB_r,
\end{eqnarray*}
 and similarly by replacing $\sigma (\beta _r)$ with
 $u(\beta _r)$, thus 
$$
 \exp \left(
 \sup\limits_{\lceil t\rceil_n \le s\le T}\int_{\lceil t\rceil_n }^s\sigma(\beta_r)dB_r+\int_{\lceil t\rceil_n }^Tu(\beta_r)dB_r-\frac{1}{2}\int_{\lceil t\rceil_n }^Tu^2(\beta_r)dr
\right)
$$
 is upper bounded by
 $$\exp
 \left(
 {\sup\limits_{t\le s\le T}\int_{t}^s(2\sigma(\beta_r)+u(\beta_r))dB_r-\inf\limits_{t\le s\le T}\int_{t}^{s}(2\sigma(\beta_r)+u(\beta_r))dB_r}
 \right),
 $$
 which is ${\P}$-integrable 
 as in \eqref{lebg3}. Therefore, by dominated convergence we find 
 \begin{eqnarray*}
   \lefteqn{\lim\limits_{s\searrow t}\E\left[e^{\sup\limits_{s\le v\le T}\int_s^v\sigma(\beta_r)dB_r+\int_s^Tu(\beta_r)dB_r-\frac{1}{2}\int_s^Tu^2(\beta_r)dr}\;\Big|\;
\beta_s=j\right]}\nonumber
\\
&=&\nonumber\lim\limits_{s\searrow t}\sum\limits_{l=1}^m
\left[e^{(s-t)Q}\right]_{j,l} \E\left[e^{\sup\limits_{s\le v\le T}\int_s^v\sigma(\beta_r)dB_r+\int_s^Tu(\beta_r)dB_r-\frac{1}{2}\int_s^Tu^2(\beta_r)dr}\;\Big|\;
  \beta_s=l\right]
\\
&=&\nonumber\lim\limits_{s\searrow t}\E\left[e^{\sup\limits_{s\le v\le T}\int_s^v\sigma(\beta_r)dB_r+\int_s^Tu(\beta_r)dB_r-\frac{1}{2}\int_s^Tu^2(\beta_r)dr}\;\Big|\;
\beta_t=j\right]\\
&=&\E\left[e^{\sup\limits_{t\le v\le T}\int_t^v\sigma(\beta_r)dB_r+\int_t^Tu(\beta_r)dB_r-\frac{1}{2}\int_t^Tu^2(\beta_r)dr}\;\Big|\;
\beta_t=j\right]
\end{eqnarray*}
and similarly,  
\begin{eqnarray}
  \label{g6}
\lefteqn{\lim\limits_{s\nearrow t}\E\left[e^{\sup\limits_{s\le v\le T}\int_s^v\sigma(\beta_r)dB_r+\int_s^Tu(\beta_r)dB_r-\frac{1}{2}\int_s^Tu^2(\beta_r)dr}\;\Big|\;
    \beta_s=j\right]
}
\\
\nonumber
&=&\E\left[e^{\sup\limits_{t\le v\le T}\int_t^v\sigma(\beta_r)dB_r+\int_t^Tu(\beta_r)dB_r-\frac{1}{2}\int_t^Tu^2(\beta_r)dr}\;\Big|\;\beta_t=j\right].
\end{eqnarray}
Combining \eqref{gx} and \eqref{g6} we conclude to Lemma \ref{gcontin}
by a classical uniform continuity argument. 
\end{Proof}
Finally, we proceed to the proof of Theorem~\ref{conver}. 
Let
\begin{equation}
\label{fjhklsjlkf} 
\Delta^n_k : = \max_{j\in {\cal M}} \sup_{a\geq 1}
| G_n(t^n_k,a,j)-G(t^n_k,a,j)|, 
\qquad k=0,1,\ldots,n,
\end{equation}
with $\Delta^n_n=0$.
By \eqref{psic0}, \eqref{th2.3} and \eqref{fjhklsjlkf}
we have
\begin{align}
  \label{deln-2}
&     \Delta^n_{k-1} \leq e^{q_{j,j}\delta_n} \Delta^n_k
    \max_{j\in {\cal M}}
    \int_0^{\infty}\int_{-\infty}^ye^{(u(j)+\sigma(j))x- u^2(j)\delta_n / 2 }
    \varphi_{\delta_n}(x,y)dxdy 
  \\
  \nonumber
  & +
  \max_{j\in {\cal M}}
  \sup_{a\geq 1}
  \sum_{i=1\atop i \not= j}^m q_{j,i}
  \int_0^{\delta_n} e^{q_{jj}r}
  \int_0^{\infty}\int_{-\infty}^ye^{(u(j)+\sigma(j))x- u^2(j)r/2}
    \\
    \nonumber
    & 
\times      | 
      G_n(t^n_k,a\vee (\sigma(j)y)-\sigma(j)x,i)
      -G( t_{k-1}^n+r,a\vee (\sigma(j)y)-\sigma(j)x,i)
    | 
    \varphi_{\delta_n} (x,y) dxdydr, 
\end{align} 
$k=1, 2, \ldots,n$, where
\begin{eqnarray}\label{gsplit}
  \nonumber
  \lefteqn{
    \! \! 
    | G_n(t^n_k,a\vee (\sigma(j)y)-\sigma(j)x,i)
      -G( t_{k-1}^n+r,a\vee (\sigma(j)y)-\sigma(j)x,i) | 
    }
  \\
  \nonumber
        & \leq & | G_n( t_k^n,a\vee (\sigma(j)y)-\sigma(j)x,i)-G( t_k^n,a\vee (\sigma(j)y)-\sigma(j)x,i) |
  \\
  \nonumber
      & & + |
      G( t_k^n,a\vee (\sigma(j)y)-\sigma(j)x,i)-G(t_{k-1}^n+r,a\vee (\sigma(j)y)-\sigma(j)x,i) |
\\  
      \nonumber
      &\leq & \Delta^n_k+\varepsilon^n_{k-1},
      \qquad
      k=1,2,\ldots,n,
      \quad a\geq 1, 
  \end{eqnarray} 
 with 
$$\varepsilon^n_k :=
\max_{i\in{\cal M}}\sup_{a\geq 1\atop t^n_k\leq s<t\leq t^n_{k+1}}|G(t,a,i)-G(s,a,i)|,
\qquad
 k=0,1,\ldots,n-1. 
$$
 Combining \eqref{deln-2} and \eqref{gsplit} yields 
 \begin{eqnarray}
\nonumber 
     \Delta^n_{k-1} & \leq & 
     \max_{j\in {\cal M}} e^{(q_{j,j}+\mu(j))\delta_n}\Delta^n_k
     \\
     \nonumber
     & & +
  (\Delta^n_k+\varepsilon^n_{k-1})
     \max_{j\in {\cal M}}
     \sum_{i=1\atop i \not= j}^m q_{j,i}
\int_0^{\delta_n} e^{q_{jj}r}
\int_0^{\infty}\int_{-\infty}^ye^{(u(j)+\sigma(j))x- u^2(j)r/2}
\varphi_{\delta_n} (x,y) dx dy dr
   \\
   \nonumber
   &= &
   e^{(q_{j,j}+\mu(j))\delta_n}\Delta^n_k +
(\Delta^n_k+\varepsilon^n_{k-1})
   \max_{j\in {\cal M}}
   \int_0^{\delta_n}\E\left[e^{(u(j)+\sigma(j))B_{\delta_n} - u^2(j)r/2+q_{jj}r}\right]dr
 \sum_{i=1\atop i \not= j}^m q_{j,i}
 \\
 \nonumber
 &= &
    e^{(q_{j,j}+\mu(j))\delta_n}\Delta^n_k +
(\Delta^n_k+\varepsilon^n_{k-1})
    \max_{j\in {\cal M}}
    \int_0^{\delta_n}
 e^{(u(j)+\sigma(j))^2\delta_n/2 - u^2(j)r/2+q_{jj}r}
 dr
 \sum_{i=1\atop i \not= j}^m q_{j,i}
 \\
\nonumber 
   & \leq & c \left(
  \Delta^n_k+\varepsilon^n_{k-1}\delta_n\right),
  \qquad
  k=1,2,\ldots,n-1,   
 \end{eqnarray}
 for some constant $c>0$ independent of $n\geq 1$, hence 
$$ 
  \Delta^n_k \leq c \left(
  \delta_n \sum_{i=k}^{n-1}\varepsilon^n_i \right),
  \qquad
  k=0,1,\ldots , n,
  $$
 and 
$$ 
 \max_{k=0,1,\ldots , n}
 \Delta^n_k =
 \max_{k=0,1,\ldots , n \atop j\in {\cal M}}
 \sup_{a\geq 1}
 | G_n(t^n_k,a,j)-G(t^n_k,a,j)|
 \leq 
 c \left(
  \max_{k=0,\ldots,n-1}\varepsilon^n_k \right)
  $$
  which tends to $0$ as $n$ tends
  to infinity by \eqref{psilim} in Lemma \ref{gcontin}.
  Consequently we have
  $$\lim\limits_{n\rightarrow\infty}
  \sup_{a\geq 1}
  |
G(\lceil t\rceil_n,a,j)-G_n(\lceil t\rceil_n,a,j) |=0$$
for any $0\leq t\le T$, $j\in {\cal M}$,
and by Lemma \ref{gcontin} it follows that
$$G(t,a,j)=\lim\limits_{n\rightarrow\infty}G(\lceil t\rceil_n,a,j)=\lim\limits_{n\rightarrow\infty}G_n(\lceil t\rceil_n,a,j),
$$
uniformly in $a\geq 1$, for all $j\in {\cal M}$ and $t\in [0,T]$.
\noindent
\hfill $\Box$

\vskip-0.4cm
\noindent
\section{Numerical results}
\label{last}
\setcounter{equation}{0}
In this section we present numerical estimates
obtained from Theorems~\ref{shutdown} and \ref{conver} for
 the boundary functions 
$$
 b_D (t,j): =\inf\{x \in [1,\infty) \ : \ (t,x,j)\in D \},
 \quad t\in [0,T], \quad j \in {\cal M},
$$
 of the stopping set $D$ defined in \eqref{setd},
 in the case of two-state Markov chains with ${\cal M} = \{ 1,2\}$.
 \\
 
 \noindent
 $(i)$ Constant drift. 
\par 

\begin{figure}[ht!]
\centering
\includegraphics[height=0.3\textwidth,width=0.8\textwidth]{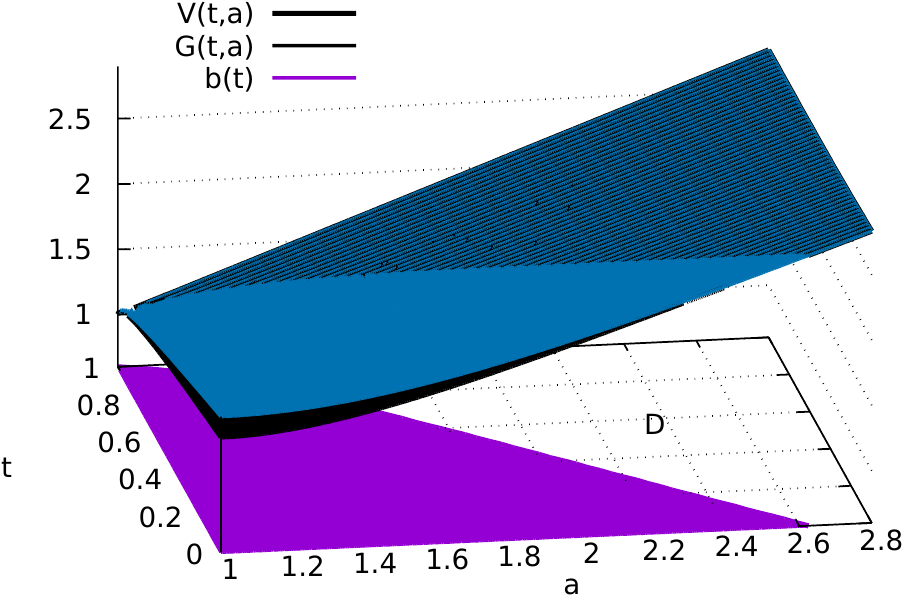}
\vskip-0.3cm
	\caption{Value functions computed from Theorems~\ref{shutdown} and \ref{conver}.}
	\label{fig0_0}
\end{figure}

\noindent 
 In the absence of regime switching, the recursive algorithm
 of Theorems~\ref{shutdown} and \ref{conver}
 is applied in Figure~\ref{fig0_0} 
 to the computation of the value functions
 $V(t,a,j)$ and $G(t,a,j)$ 
 with $T=1$, $\sigma =0.5$, $\mu =0.2$, $n=50$, and $\delta_n = T/n = 0.01$. 
 
 Figure~\ref{fig0_0} allows us in particular to visualize
 the stopping set $D$ defined in \eqref{setd}
 and the continuation set
 $C = \big\{
 (t,a)\in[0,T] \times [1,\infty) \ : \ V(t,a) < G(t,a)
 \big\}$. 
 \\

\begin{figure}[ht!]
\centering
\includegraphics[height=0.3\textwidth,width=0.8\textwidth]{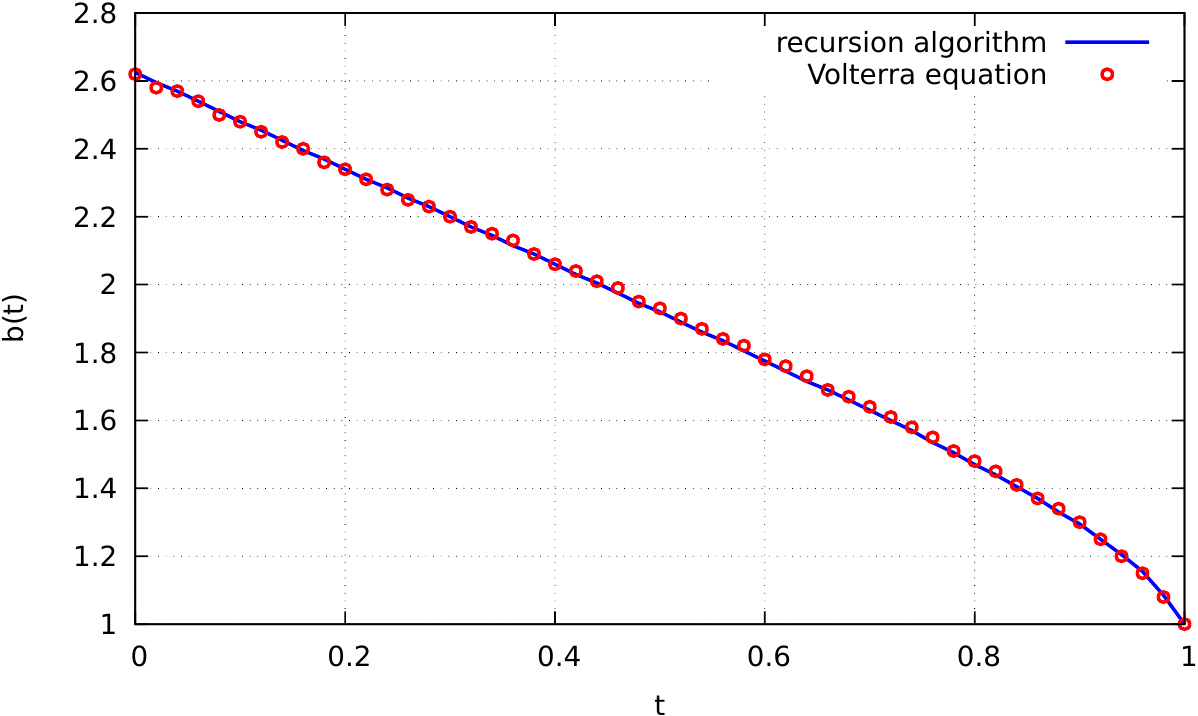}
\vskip-0.2cm
	\caption{Boundary function computed from Theorems~\ref{shutdown} and \ref{conver} {\em vs} \eqref{ebdry}.}
	\label{fig0}
\end{figure}

 \noindent
In Figure~\ref{fig0} the recursive method is compared to 
the solution of the Volterra integral equation
\eqref{ebdry} by dichotomy for the computation
of the boundary function $b(t)$. 

\pagebreak[4]

\noindent 
As shown in Figure~\ref{fig0}, the recursive and Volterra equation
methods yield similar levels of precision.
However, increasing the number $n$ of time steps will make the Volterra
equation method perform slower relative to the recursion method,
due to the quadratic complexity of the former and to the linear complexity of the latter. 
\\
 
\noindent
 $(ii)$ Drifts with switching signs. 
 \par
 
\begin{figure}[ht!]
  \hskip-0.2cm
  \begin{subfigure}{.4\textwidth}
  \centering
  \includegraphics[width=7.8cm,height=4.5cm]{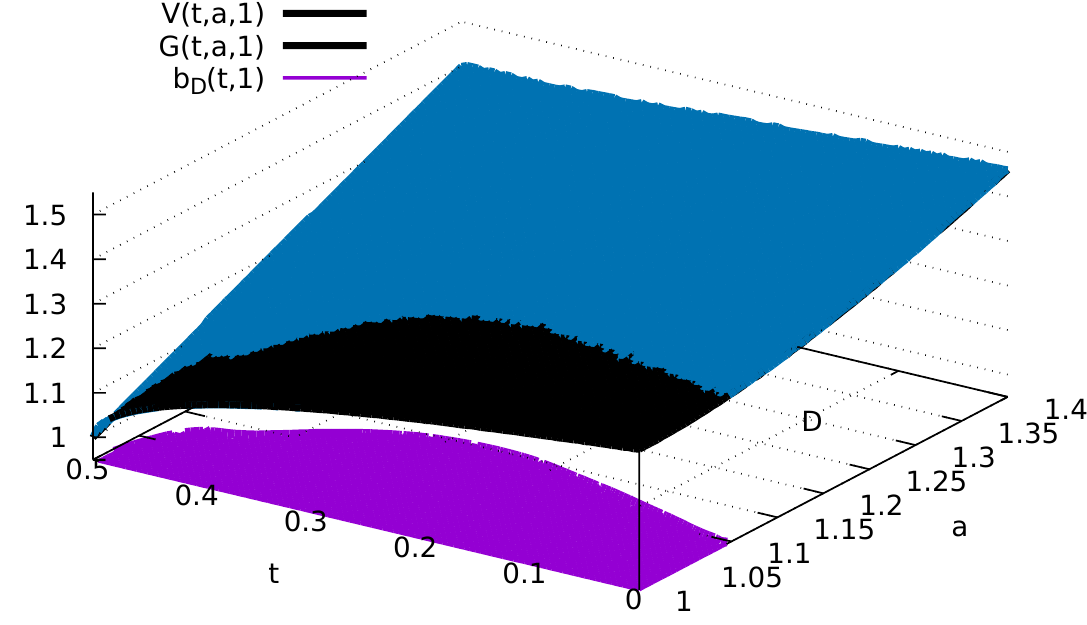}
\end{subfigure}
\hskip2cm
\begin{subfigure}{.4\textwidth}
  \centering
  \includegraphics[width=7.8cm,height=4.5cm]{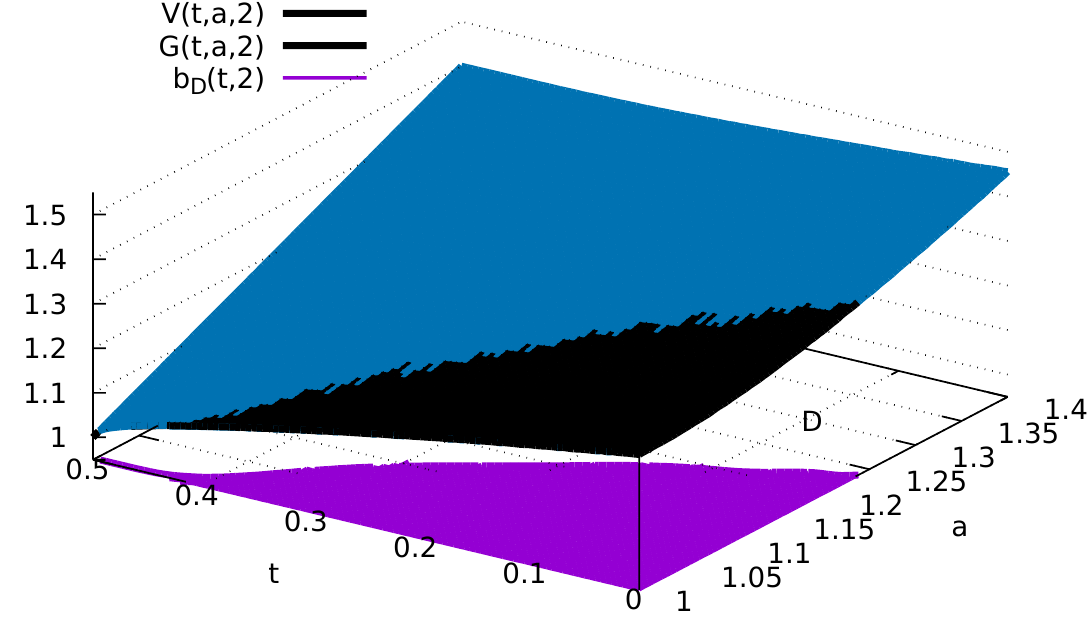}
\end{subfigure}
\vskip-0.1cm
\caption{Value functions under drifts of mixed signs.}
\label{fig2-3}
\end{figure}

 \noindent 
 Figure~\ref{fig2-3} presents the graphs of the value
 functions obtained from
 the recursive algorithm of Theorems~\ref{shutdown} and \ref{conver}
 with $\mu (1) =0.2$, $\mu (2) =-0.2$, $\sigma (1) =0.5$,
 $\sigma (2) =0.3$, $T=0.5$, $n=100$, $\delta_n = T/n = 0.05$, 
 and 
\begin{displaymath}
\mathbf{Q} =
\left[ \begin{array}{cc}
-2.5 & 2.5
\\
2 & -2 \nonumber
\end{array}\right].
\end{displaymath}
  
\noindent

Figure~\ref{fig2-3} also allows us to visualize
the stopping set $D$ and the continuation set 
\begin{equation*}
 C = \big\{
 (t,a,j)\in[0,T] \times [1,\infty) \times {\cal M} \ : \ V(t,a,j) < G(t,a,j)
   \big\}.
\end{equation*} 
The numerical instabilities observed are due to the necessity
to check the equality $V(t,a,j) = G(t,a,j)$ when
$V(t,a,j)$ and $G(t,a,j)$ are very close to each other.
We observe that the corresponding boundary 
function $t \mapsto b_D(t,1)$ starting from state $1$ is not monotone.
Precisely, when time $t$ is close to $0$ 
it is better to exercise early because one may 
switch to state $2$ after the average time 
$1/q_{1,1} = 0.4$, 
in which case the drift takes the negative value $\mu (2)=-0.2$.
On the other hand, when $t$ increases up to $0.3$ the function
$t \mapsto b_D(t,1)$ tends to increase as it makes more sense to wait since 
we may stay at state $1$ with $\mu(1) =0.2$ for the remaining average time
$1/q_{2,2}=0.5$, which is higher than the remaining time $T-t$.

\begin{figure}[ht!]
\centering
\includegraphics[height=0.3\textwidth,width=0.8\textwidth]{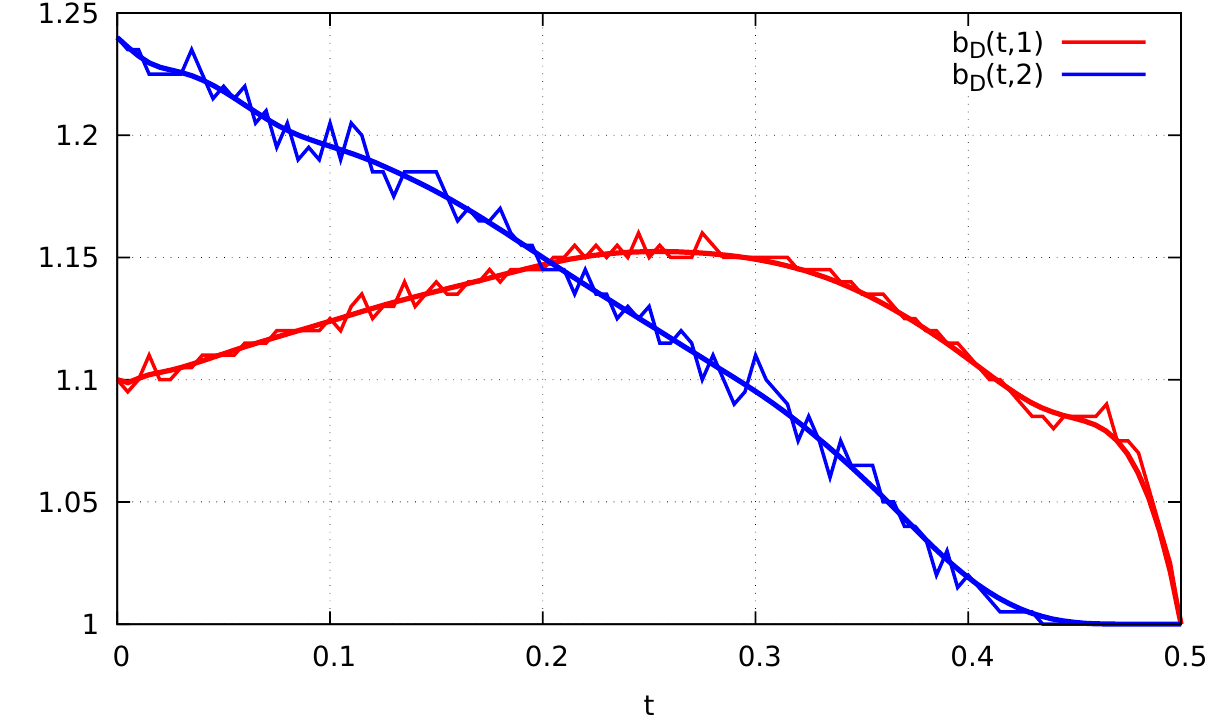}
\vskip-0.2cm
	\caption{Boundary functions under drifts of mixed signs.}
	\label{fig3}
\end{figure}

The boundary functions are plotted in Figure~\ref{fig3} 
with spline smoothing. 
Starting from state $2$ we observe the usual decreasing boundary
$t \mapsto b_D(t,2)$,
which here becomes close to $0$ 
before time $T$, since in this case
we should exercise immediately as the
average time $1/q_{2,2}=0.5$ to switch to state $1$
exceeds the remaining time $T-t$ until maturity.
\\ 

Until time $0.2$ we should exercise
immediately when switching from state $2$ to state $1$
at a time $t$ such that $b(t,1) < \hat{Y}_{0,t} / Y_t = a < b_D(t,2)$,
while after time $0.2$ the strategy is the opposite if
$b(t,2) < \hat{Y}_{0,t} / Y_t = a < b_D(t,1)$.
 
\footnotesize 

\def\cprime{$'$} \def\polhk#1{\setbox0=\hbox{#1}{\ooalign{\hidewidth
  \lower1.5ex\hbox{`}\hidewidth\crcr\unhbox0}}}
  \def\polhk#1{\setbox0=\hbox{#1}{\ooalign{\hidewidth
  \lower1.5ex\hbox{`}\hidewidth\crcr\unhbox0}}} \def\cprime{$'$}

\end{document}